\documentclass{amsart}
\usepackage{amsmath,amsthm}
\usepackage{amsfonts,amssymb}

\newtheorem{thm}{Theorem}[section]

\newtheorem{lem}[thm]{Lemma}
\newtheorem{prop}[thm]{Proposition}
\newtheorem{exam}[thm]{Example}

\theoremstyle{remark}

 \def\ab{{\mathbf a}}

 \def\tb{{\mathbf t}}

 \def\NN{{\mathbb N}}

 \def\RR{{\mathbb R}}
 \def\SS{{\mathbb S}}
 
        \def\sspan{\operatorname{span}}
        
        \def\dim{\operatorname{dim}}

\newcommand{\wt}{\widetilde}

\begin{document}
 
\title{Polynomial Interpolation on the Unit Sphere II}
\author{ Wolfgang zu Castell, Noem\'{\i} La\'{\i}n Fern\'{a}ndez and Yuan Xu} 
\address{Institute of Biomathematics and Biometry\\
      GSF - National Research Center for Environment and Health\\
      85764 Neuherberg, Germany} \email{castell@gaf.de}
\address{Center for Mathematical Sciences \\Munich University of Technology\\Boltzmannstr. 3, 85747 Munich, Germany}\email{fernande@ma.tum.de}
\address{Department of Mathematics\\ University of Oregon\\
    Eugene, Oregon 97403-1222.}\email{yuan@math.uoregon.edu}

\date{\today}
\thanks{The second author was supported by the Graduate Program 
\emph{ Applied Algorithmic Mathematics} of the Munich University of 
Technology. The work of the third author was supported in part by the 
National Science Foundation under Grant DMS-0201669}
                          
\begin{abstract}
The problem of interpolation at $(n+1)^2$ points on the unit sphere 
$\mathbb{S}^2$ by spherical polynomials of degree at most $n$ is proved 
to have a unique solution for several sets of points. The points are 
located on a number of circles on the sphere with even number of points on 
each circle. The proof is based on a method of factorization of polynomials. 
\end{abstract}

\maketitle                      
 
\section{Introduction}
\setcounter{equation}{0}

Let $\SS^2 = \{x: \|x\| =1\}$ denote the unit sphere of $\RR^3$, where 
$\|x\|^2 = x_1^2+x_2^2+x_3^2$. Let $\Pi_n(\SS^2)$ denote the space of 
spherical polynomials of degree $n$, which is the restriction  
of  polynomials of degree  $n$ in three variables to $\SS^2$. It is known that 
$$
    \dim \Pi_n(\SS^2) = (n+1)^2,   \qquad n \ge 0.
$$ 
The problem of interpolation on the unit sphere by polynomials is as
follows: 

\medskip\noindent 
{\bf Problem 1.} Let $X = \{\ab_i: 1 \le i \le (n+1)^2\}$ be a set of pairwise distinct 
points on $\SS^2$. Find conditions on $X$ such that there is a unique 
polynomial $T \in \Pi_n(\SS^2)$ satisfying 
$$
    T(\ab_i) = f_i, \qquad \ab_i \in X, \quad 1 \le i \le (n+1)^2,
$$
where $\{f_i\}$ is an arbitrary set of data.

\medskip
If there is a {\it unique} solution to the interpolation problem, we say that 
the problem is {\it poised} and that $X$ solves Problem 1. This problem has 
been studied recently in \cite{F1,F,FP,GL,SW,X1,X2}. 

Although almost all choices of $X$ will solve Problem 1, it is difficult to 
know whether a given set $X$ will work since  computing the determinant 
of the  interpolation matrix is difficult. In \cite{X1} a large family of sets of 
interpolation points is given explicitly, each set solving  Problem 1. Let us briefly describe this construction. The 
$(n+1)^2$ points lie on $n+1$ distinct latitudes (parallel circles on 
$\SS^2$), and each latitude contains an odd number of equidistant points. 
The number of points needs not to be the same on each latitude and there is no 
restriction on the position of the latitudes. For the simplest case   
$n = 2m$, the set of $(2m+1)^2$ points  lie on $2m+1$ latitudes, each of them containing  $2m+1$ 
equally spaced nodes. In \cite{F}, another family of points 
that solves Problem 1 was found, for which $n = 2m-1$. There the points lie on $2m$ 
latitudes and each latitude has an even number of $2m$ equally spaced points.
In this case, the $2m$ latitudes are divided into two groups; the equidistant
points on one half of the latitudes need to differ by a rotation from the 
points on the other half of the latitudes. 
While the proof in \cite{F} is based on the analysis of  the determinants of the interpolation 
matrix, the proof in \cite{X1} uses a factorization method which
avoids the determinants. Furthermore, the factorization method provides many 
more sets of points leading to poised problems. A key observation in \cite{X1} is 
that the use of equidistant points allows us to reduce the problem on the sphere to a special trigonometric interpolation problem. 

The purpose of this paper is to show that the factorization method also 
works in the setting of an even number of points on each latitude. Again, the 
use of equidistant points reduces the problem to an interpolation problem 
of one variable. However, the new interpolation problem is different from 
the one with an  odd number of points on each latitude and has to be solved using a 
completely different method. In comparison to \cite{F},  the factorization method allows to  obtain
more sets of points that solve Problem 1. 

For the background of polynomial interpolation in general, we refer to 
the survey article \cite{GS} and the references therein, even though
interpolation on the sphere is not discussed there. Let us also mention 
that the factorization method is closely related to the method used for
polynomial interpolation on the unit disk in \cite{BX1,BX2,HS}. Apart from
a result in \cite{GL}, which is a simple consequence of Bezout's theorem, the family
of points found in \cite{F,X1}, and those stated below appear
to be the only ones that are given explicitly for all $n$. 

The paper is organized as follows. The factorization  method  is studied in 
Section 2 and its application to polynomial interpolation on $\SS^2$ is given
in Section 3. 

\section{Factorization of polynomials}
\setcounter{equation}{0}

\subsection{Polynomial representation}
For fixed $a\in (-1,1)$, let $\mathbb{S}^2(a):=\lbrace (x,y,z)\in
\mathbb{S}^2:\,z=a\rbrace$ denote the circle on $\mathbb{S}^2$ resulting 
from the intersection of $\mathbb{S}^2$ with the plane $z=a$. This set is 
called \emph{latitude} at $z=a$. 

On the unit sphere $\SS^2$ it is more convenient to work with spherical 
coordinates, 
\begin{align*}
  x = \sin \theta \sin \phi, \quad y = \sin \theta \cos \phi, \quad 
  z = \cos \theta, \qquad 0 \le \phi < 2 \pi, \quad 0 \le \theta \le \pi.
\end{align*} 
For a polynomial $T_n\in\Pi_n(\mathbb{S}^2)$, 
we introduce the notation $\wt T_n$ defined by 
\[
\wt T_n(\theta,\phi)=T_n(\sin\theta\,\cos\phi,\;\sin\theta\,\sin\phi,\;
 \cos\theta),\quad 0\le \phi< 2\pi,\quad 0\le \theta\le\pi.
\]
If $X\!=\!\lbrace (x_i,y_i,z_i):\,1\le i\le M \rbrace$ is a set of points on 
$\mathbb{S}^2$, we also use the notation 
$\wt X\!=\!\lbrace (\theta_i,\phi_i):\,1\le i\le M \rbrace$ for the 
corresponding set of spherical coordinates.

It has been shown in Section 2 of \cite{X1} that the polynomial $\wt T_n$ 
can be written as
\begin{align}\label{eq:polynomial}
& \wt T_n(\theta,\phi)= a_0(\cos\theta) \\
& \qquad\qquad +\sum\limits_{k=1}^{n}
  \left[a_k(\cos\theta)\,(\sin\theta)^k\cos k\phi
         +b_k(\cos\theta)(\sin\theta)^k\sin k\phi \right], \notag
\end{align}
where $a_k(\cdot)$ and $b_k(\cdot)$ are polynomials of degree $n\!-\!k$ in 
one variable. Note that for any fixed $\theta$, the polynomial 
$\wt T_n(\theta,\cdot)$ is a trigonometric polynomial of degree $n$.

Below we will consider interpolation problems based on points that are 
equidistantly distributed  on an even number of latitudes, each of them  containing an even number 
of nodes. To describe these points, it is convenient to introduce the 
following notation: 
\[
\Theta_{\alpha,s}:=\left\lbrace \phi_j^{\alpha}: 
  \phi_j^{\alpha}=\frac{(2j+\alpha)\pi}{2s},\ j=0,1,\dots,2s-1 \right\rbrace,
\] 
where $s\in\mathbb{N}$ and $ \alpha\in [0,2)$. These points can be considered as equidistant points on the 
unit circle using the mapping 
$\phi\mapsto e^{i\phi}$. The parameter $\alpha$ indicates that the points 
are defined up to a rotation by an angle of $\alpha\pi/2s$. 

\begin{lem}  \label{lem:2.1} 
Let $n=2m\!-\!1$ and $\alpha\in [0,2)$. For $\phi\in \Theta_{\alpha,m}$, 
\begin{align} \label{eq:representation}
& \wt T_n(\theta,\phi)=  a_0(\cos\theta) \\
& \qquad + \sum\limits_{k=1}^{m-1}\left[ (a_k(\cos\theta)(\sin\theta)^k\!+
   \!u_{2m-k}(\cos\theta)(\sin\theta)^{2m-k})\cos k\phi \right . \notag\\
& \qquad\qquad \quad \left . +(b_k(\cos\theta)(\sin\theta)^k\!+ 
 \!v_{2m-k}(\cos\theta)(\sin\theta)^{2m-k})\sin k\phi \right]\notag \\
& \qquad + 
 \left(a_m(\cos\theta)\cos\frac{\alpha\pi}{2}-b_m(\cos\theta)
\sin\frac{\alpha\pi}{2}\right)(\sin\theta)^m 
 \cos\left(m\phi - \frac{\alpha\pi}{2}\right), \notag 
\end{align}
where, for $k=1,\dots,m-1$,
\begin{align*}
u_{2m-k}(t) & =
   a_{2m-k}(t)\cos\alpha\pi+b_{2m-k}(t)\sin\alpha\pi,\\
v_{2m-k}(t) & =a_{2m-k}(t)\sin\alpha\pi-b_{2m-k}(t)\cos\alpha\pi
\end{align*}
are polynomials of degree $k-1$. 
\end{lem}

\begin{proof} 
We split the sum in \eqref{eq:polynomial} into two sums, one over 
$1 \le k \le {m-1}$ and the other over $m \le k \le 2m-1$. In the second
sum we change the summation index $k \mapsto 2m\!-\!k$ and use the
elementary relations 
\begin{align*}
\cos (2m-k)\phi = & \cos ((2j+\alpha)\pi-k\phi)=\cos\alpha\pi\cos k\phi+
  \sin\alpha\pi\sin k\phi,\\
\sin (2m-k)\phi= & \sin ((2j+\alpha)\pi-k\phi)=\sin\alpha\pi\cos k\phi-
  \cos\alpha\pi\sin k\phi,
\end{align*}
which holds for $\phi\in\Theta_{\alpha,m}$. Combining the two sums, we obtain 
\begin{align*} 
\wt T_n(\theta,\phi)= & \ a_0(\cos\theta) 
 + \sum\limits_{k=1}^{m-1}\left[ (a_k(\cos\theta)(\sin\theta)^k\!+
   \!u_{2m-k}(\cos\theta)(\sin\theta)^{2m-k})\cos k\phi \right . \notag\\
& \qquad\qquad \qquad \left . +(b_k(\cos\theta)(\sin\theta)^k\!+ 
 \!v_{2m-k}(\cos\theta)(\sin\theta)^{2m-k})\sin k\phi \right]\notag \\
& \qquad + 
(a_m(\cos\theta)\cos (\alpha\pi-m\phi)+
    b_{m}(\cos\theta)\sin (\alpha\pi-m\phi))(\sin\theta)^m. \notag
\end{align*}
Using the addition formula for the cosine and the sine function
\begin{align*}
\cos(m\phi-\alpha\pi) = & 
\cos\left (m\phi-\frac{\alpha\pi}{2}\right)\cos\frac{\alpha\pi}{2}+
 \sin\left (m\phi-\frac{\alpha\pi}{2}\right)\sin\frac{\alpha\pi}{2},\\
\sin(m\phi-\alpha\pi) = & \sin\left (m\phi-\frac{\alpha\pi}{2}\right)
  \cos\frac{\alpha\pi}{2}-\cos\left (m\phi-\frac{\alpha\pi}{2}\right)
  \sin\frac{\alpha\pi}{2},
\end{align*}
the $a_m$ and $b_m$ terms of the above expression of $\wt T_{n}$ 
can be rewritten as 
\begin{align*} 
& \left(a_m(\cos\theta)\cos\frac{\alpha\pi}{2}-b_m(\cos\theta)
\sin\frac{\alpha\pi}{2}\right)(\sin\theta)^m 
 \cos\left(m\phi - \frac{\alpha\pi}{2}\right)\\
+ & \left (a_m(\cos\theta)\sin\frac{\alpha\pi}{2}
  +b_m(\cos\theta)\cos\frac{\alpha\pi}{2}\right)(\sin\theta)^m
  \sin\left(m\phi- \frac{\alpha\pi}{2}\right).
\end{align*}
Consequently,  formula \eqref{eq:representation} follows from the
fact that $\phi \in \Theta_{\alpha,m}$ satisfies $\sin (m \phi- \alpha\pi/2) 
=0$.
\end{proof}

\begin{lem} \label{lem:2.2} 
Let $\wt T_{2m-1}$ be given as in \eqref{eq:representation} and  $\theta\in (0,\pi)$. If 
$\wt T_{2m-1}(\theta, \phi) = 0$ for $\phi \in \Theta_{\alpha,m}$, then
$a_0(\cos\theta)=0$, and  
\begin{align} \label{eq:2.3}
\begin{split}
& a_k(\cos\theta)+(\sin\theta)^{2m-2k}(a_{2m-k}(\cos\theta)\cos\alpha\pi
 +b_{2m-k}(\cos\theta)\sin\alpha\pi)= 0,\\
& b_k(\cos\theta)+(\sin\theta)^{2m-2k}(a_{2m-k}(\cos\theta)\sin\alpha\pi-
   b_{2m-k}(\cos\theta)\cos\alpha\pi) = 0,
\end{split}
\end{align}
for $1 \le k \le m-1$. Furthermore,
\begin{equation}\label{eq:2.4}
  a_m(\cos\theta)\cos\frac{\alpha\pi}{2}
    -b_m(\cos\theta)\sin\frac{\alpha\pi}{2} =0.
\end{equation} 
\end{lem}

\begin{proof}
The proof uses the following fact. Interpolation on the $2m$ zeros of 
the function 
$\sin\left(m\phi -\frac{\alpha\pi}{2}\right)$ inside $[0, 2 \pi)$ by a 
trigonometric polynomial of the form 
$$
S_{m-1}(\phi) = a_0 + \sum_{k=1}^{m-1} \left(a_k \cos k \theta + 
  b_k \sin k \theta \right) + a_n \cos (m\phi-\alpha\pi/2)
$$
is unique (see \cite[Vol. II]{Z}). Observe that the points in 
$\Theta_{\alpha,m}$ are exactly zeros of $\sin\left(m\phi-\alpha\pi/2\right)$
and $\cos (m\phi-\alpha\pi/2) = \pm 1$ for $\phi \in \Theta_{\alpha,m}$.

By Lemma \ref{lem:2.1}, $\wt T_{2m-1}$ takes the form of 
\eqref{eq:representation}. The assumption that $\wt T_{2m-1}(\theta, \phi)
\!\!=\!\!0$ implies that the coefficients of $\wt T_{2m-1}(\theta, \cdot)$ are all 
zero. This leads to $a_0(\cos \theta) =0$, 
\begin{align*}
& a_k(\cos\theta)(\sin\theta)^k+u_{2m-k}(\cos\theta)(\sin\theta)^{2m-k}=0,\\
& b_k(\cos\theta)(\sin\theta)^k+v_{2m-k}(\cos\theta)(\sin\theta)^{2m-k}=0,
\end{align*}
for $k=1,\dots,m-1$, which become the equations in \eqref{eq:2.3} upon 
multiplying by $(\sin \theta)^{-k}$, and 
$$
\left(a_m(\cos\theta)\cos\frac{\alpha\pi}{2}-b_m(\cos\theta)
\sin\frac{\alpha\pi}{2}\right)(\sin\theta)^m=0, 
$$
which gives \eqref{eq:2.4}. Note that  $\sin\theta\neq 0$ as $\theta\in (0,\pi)$. 
\end{proof}

To proceed from here, we want to choose $2m$ distinct  $\theta_i$ 
such that whenever the polynomials in \eqref{eq:2.3} and \eqref{eq:2.4} 
vanish on these $2m$ points, they will be identically zero. To this end, 
however, we need to impose an additional symmetry. We choose 
$\theta_i$ to satisfy 
\begin{equation} \label{eq:theta}
  \theta_{2m+1-i}=\pi-\theta_i, \qquad 
     \theta_i \in (0, \pi), \quad 1 \le i \le m. 
\end{equation}
In other words, we choose the latitudes to be symmetric with respect to
the equator. 

If $p(t)$ is a polynomial of degree $n$, we denote by $p^{{\rm even}}$
and $p^{{\rm odd}}$ the even and the odd part of $p$, respectively. To be precise, if $p (t) = \sum_{j=0}^N a_j t^j$, then
$$
p^{{\rm even}}(t) = \sum_{0 \le 2j \le N} a_{2j} t^{2j} 
  \quad \hbox{and} \quad  
p^{{\rm odd}}(t) = \sum_{1 \le 2 j-1\le N} a_{2j-1} t^{2j-1}.
$$

\begin{lem} \label{lem:2.3}
Let $\wt T_{2m-1}$ be given as in \eqref{eq:representation}. If for
some $\theta$ in $(0, \pi)$  
$$
\wt T_{2m-1}(\theta, \phi) = 0, \quad \hbox{$\phi \in \Theta_{0,m}$ and}
\quad \wt T_{2m-1}(\pi - \theta, \phi) = 0, \quad \phi \in \Theta_{1,m}, 
$$
then, setting $t= \cos \theta$, we have $a_0(t)=0$ and, for $1 \le k \le m-1$, 
\begin{align} \label{eq:2.6}
\begin{split}
& p_{2m-k-1}^{{\rm even}}(t)+ q_{k-1}^{{\rm odd}}(t)
    (1-t^2)^{m-k} = 0, \\
& p_{2m-1-k}^{{\rm odd}}(t)+q_{k-1}^{{\rm even}}(t)(1-t^2)^{m-k}=0, 
\end{split}
\end{align}
where either $p_{2m-k-1}(t) = a_k(t)$ and $q_{k-1}(t) = a_{2m-k}(t)$, or 
$p_{2m-k-1}(t) = b_k(t)$ and $q_{k-1}(t) = - b_{2m-k}(t)$. Furthermore,
$ a_m(t) = b_m (t) =0$.
\end{lem}

\begin{proof}
The assumption allows us to use the previous lemma. Since $\alpha =0$ or 
$\alpha =1$, the fact that $a_m(t) = b_m(t) = 0$ follows immediately from 
\eqref{eq:2.4}. For $1 \le k \le m-1$ and $\alpha=0$, the  equations \eqref{eq:2.3}  become
\begin{align}  \label{eq:2.7}
\begin{split}
& a_k(\cos\theta)+(\sin\theta)^{2m-2k} a_{2m-k}(\cos\theta)= 0,\\
& b_k(\cos\theta) - (\sin\theta)^{2m-2k} b_{2m-k}(\cos\theta) = 0.
\end{split}
\end{align}
For  $1 \le k \le m-1$, $\alpha=1$ and $\theta$ replaced by
$\pi - \theta$, the equations \eqref{eq:2.3} take the form 
\begin{align}  \label{eq:2.8}
\begin{split}
& a_k(-\cos\theta) - (\sin\theta)^{2m-2k} a_{2m-k}(-\cos\theta)= 0,\\
& b_k(-\cos\theta) + (\sin\theta)^{2m-2k} b_{2m-k}(-\cos\theta) = 0.
\end{split}
\end{align}
Since $ p (t) + p(-t) = 2 p^{{\rm even}}(t)$ and $ p (t) - p(-t) = 
2 p^{{\rm odd}}(t)$,  combining    equations  \eqref{eq:2.7} and 
\eqref{eq:2.8} proves the result. 
\end{proof}

We use the notation $p_{2m-k-1}$ and $q_{k-1}$ since they are polynomials 
of degree $2m-k-1$ and $k-1$, respectively. In the following we will 
work with the equations in \eqref{eq:2.6}. If $p$ is an even polynomial, 
 it can be written as $p(t) = p^*(t^2)$; if $p$ is an odd polynomial, 
 it can be written as $p(t) = t p^*(t^2)$. Thus, in place of 
\eqref{eq:2.6} we will need to consider polynomials  of the form $p(t) + t q(t) (1-t^2)^r$ and $t p(t) + q(t) (1-t^2)^r$. We will need to 
study the possibility of interpolation by such polynomials. This is 
discussed in the following subsection. 

\subsection{Chebyshev systems}

A family of functions $\{\phi_1,\ldots,\phi_r\}$ is called a Chebyshev
system on a set $E \subseteq \RR$, if every linear combination from 
 the $\sspan \{\phi_1,\ldots,\phi_r\}$ has at most $r$ zeros in $E$; in 
other words,  interpolation on $r$ points by functions in 
 the $\sspan \{\phi_1,\ldots,\phi_r\}$ has a unique solution. In this subsection
we prove that the families of functions  in \eqref{eq:2.6} are 
Chebyshev systems on $(0,1)$.

\begin{prop}\label{prop:1}
Let $r$ and $s$ be two nonnegative integers such that $r >s >0$. For 
$\epsilon = 0$ or $1$, let 
$$
g(t) = p_r(t^2) + t^{\pm 1} (1-t^2)^{r-s} q_{s-1+\epsilon}(t^2),
$$ 
where $p_r$ and $q_{s-1+ \epsilon}$ are polynomials of degree $r$ and 
$s-1+\epsilon$, respectively. 
If $g$ vanishes on $r+s+1+\epsilon$ distinct points in $(0,1)$, then 
$g(t) \equiv 0$. 
\end{prop}

\begin{proof}
We first prove the case that the power of $t^{\pm 1}$ in $g(t)$ is taken as $t$ and  $\epsilon =0$. The cases $t^{-1}$ or $\epsilon =1$ are 
similar; in fact, the proof for the case $\epsilon =1$ is identical, 
and only minor changes (merely the numbers $b_k^*$ below will change)
are needed for the case that $t^{\pm 1}$ is taken as $t^{-1}$.  

Changing variables $t \mapsto t^2$ shows that we need to prove that if 
$$
  g^*(t) = p_r(t) + \sqrt{t} (1-t)^{r-s} q_{s-1}(t) 
$$ 
vanishes on $r+s+1$ distinct points in $[0,1]$, then $g^*(t) \equiv 0$. Let 
$$
h(t):= t^{r-\frac12} \frac{d^{r+1}}{dt^{r+1}} g^*(t) = 
   t^{r-\frac12}\frac{d^{r+1}}{dt^{r+1}} \left[\sqrt{t} (1-t)^{r-s} 
     q_{s-1}(t)\right]. 
$$  
Using Rolle's  theorem repeatedly, we see that it suffices to prove that 
if  $h(t)$ vanishes on $s$ distinct points in $(0,1)$, then 
$q_{s-1} (t) \equiv 0$. 

Since $q_{s-1}$ is a polynomial, we can write it as 
$$
  q_{s-1} (t) = b_0 + b_1 (1-t) + \ldots + b_{s-1} (1-t)^{s-1}.
$$
Using the Leibnitz rule repeatedly, we have 
$$
 \frac{d^{r+1}}{dt^{r+1}} t^{k+j+\frac12} = \frac{t^{-r-\frac12}}{2^{r+1}} 
   \prod_{i=0}^{k+j} (2 i +1) (-1)^{r-k-j} 
   \prod_{i=1}^{r-k-j} (2 i -1)  t^{k+j}. 
$$
In the following, we will use the convention that $\prod_{i=a}^{b} =1$  
whenever $b < a$. This leads to 
\begin{align*}
h(t) = \frac{d^{r+1}}{dt^{r+1}}  \sum_{k=0}^{s-1} b_k 
         \sum_{j=0}^{r-s} (-1)^j \binom{r-j}{j} t^{k+j+\frac12} 
     = \sum_{k=0}^{s-1} b_k^* h_k (x),
\end{align*}
where 
$$
 b_k^* = b_k \frac{(-1)^{r-k}}{2^{r+1}} \prod_{i=0}^{k-1} (2 i +1) 
      \prod_{i=1}^{s-k-1} (2 i-1) 
 \quad\hbox{and}\quad  h_k(t) = \sum_{j=0}^{r-s} a_{k,j} t^{j+k}, 
$$
in which the coefficients $a_{k,j}$ are given by
\begin{equation}\label{eq:akj}
       a_{k,j} = \binom{r-s}{j} \prod_{i=0}^{j} (2 k+2i +1)   
              \prod_{i=j}^{r-s} (2 (r-k-i)-1). 
\end{equation}
We note that all coefficients $a_{k,j}$ are positive numbers. The
polynomial $h_k$ is of degree $r-s+k$. In order
to prove the proposition, we need to show that the set $\{h_0, h_1, 
\ldots, h_{s-1}\}$ forms a Chebyshev system on $(0,1)$. In other words,
we need to prove that the matrix $(h_j(t_k))_{j,k=0}^{s-1}$ is invertible
for any set of distinct points in $(0,1)$. 

Let $\tb = \{t_1, t_2, \ldots,t_s\}$ be a given set of distinct numbers 
in $(0,1)$. For a given set of nonnegative integers $\lambda = 
 \{j_0, j_1, \ldots, j_{s-1} \}$, we introduce the notation 
$$
   V(\lambda;\tb) = \det \left[ \begin{matrix} 
       t_1^{j_0}& t_2^{j_0}&  \ldots &  t_s^{j_0} \\ 
       t_1^{j_1}& t_2^{j_1}&  \ldots &  t_s^{j_1} \\ 
       \vdots& \vdots &  \ldots &  \vdots \\ 
       t_1^{j_{s-1}}& t_2^{j_{s-1}}&  \ldots &  t_s^{j_{s-1}} \\ 
       \end{matrix} \right].
$$
In the case of $\lambda = \{s-1,s-2,\ldots,0\}$, we denote the determinant
by $V_s(\tb)$, which is the Vandermonde determinant 
$$ 
V_s(\tb) = \det [t_{k+1}^j]_{k,j = 0}^{s-1} = 
   \prod_{1 \le i < j \le s} (t_j - t_i).  
$$    
For a given set of nonnegative integers $\lambda$, we further introduce
the notation
$$
s_\lambda(\tb) = s_{j_0,j_1,\ldots,j_{s-1}}(\tb) = 
     \frac{V(\lambda;\tb)}{V_s(\tb)}. 
$$
Note that $s_\lambda$ is a symmetric polynomial in $\tb$ and $s_\lambda$ 
is zero if $j_0, j_1,\ldots,j_{s-1}$ are not pairwise distinct. If $\mu = (\mu_0,
\mu_1\ldots,\mu_{s-1})$ is a partition, that is, $\mu_0 \ge \mu_1 \ge \ldots 
\ge \mu_{s-1} \ge 0$, $\mu_i \in \NN$, and $j_i = \mu_i + n-i+1$ for
$0 \le j \le s-1$, then $s_\lambda$ is called a Schur polynomial, cf. \cite{M} for details. It is known that Schur polynomials can be written as a linear combination 
of monomial symmetric polynomials and the coefficients in the linear 
combination are all positive (called Kosta numbers). For our purpose, it 
is enough to note that the Schur polynomials are positive when $t_l > 0$ 
for all $1 \le l \le s$. In particular, it follows that if 
$j_0 < j_1 < \ldots < j_{s-1}$, then $s_\lambda (\tb)$ is positive when
$t_l > 0$ for all $1 \le l \le s$. 

Using the definition of the polynomials $h_j$ we can write
\begin{align} \label{eq:sum}
& \frac{1}{V_s(\tb)} \det (h_j(t_{k+1}))_{j,k=0}^{s-1} \\
& \quad =  \sum_{j_0=0}^{r-s} \sum_{j_1=1}^{r-s+1} \cdots  
  \sum_{j_{s-1}={s-1}}^{r-1} 
  a_{0,j_0} a_{1,j_1-1} \ldots a_{s-1,j_{s-1}-(s-1)}    
   s_{j_0,j_1,\ldots, j_{s-1}}(\tb). \notag   
\end{align}
For $\lambda = \{j_0,j_1,\ldots,j_{s-1}\}$ we will also denote the 
coefficient of $s_\lambda$ in the above sum by $A_\lambda$; that is 
$$
  A_\lambda = a_{0,j_0} a_{1,j_1-1} \ldots a_{s-1,j_{s-1}-(s-1)}.    
$$
Recall that $a_{j,k}$ are all positive numbers. If $\lambda = 
\{j_0,j_1,\ldots,j_{s-1}\}$ is not a partition, then a 
proper permutation of $j_0,j_1,\ldots,j_{s-1}$ will be. The determinant 
changes sign when two rows are exchanged, so that $s_\lambda$ is positive 
if the permutation is even and it is negative if the permutation is odd. 
Every permutation can be factored into a number of transpositions. A
transposition means exchanging two elements. If $s_\lambda$
is negative, there is a transposition of $\lambda$, call it $\lambda'$, 
such that $s_{\lambda'} (\tb) =  - s_\lambda(\tb) > 0$. 

Let $s_\lambda(\tb)$ be negative and assume that $j_p$ and $j_q$ are a pair 
in $\lambda = \{j_0, j_1,\ldots,j_{s-1}\}$ such that $p < q$ but $j_p > j_q$. 
Considering  the summation indices in \eqref{eq:sum}, we must 
have $j_p > j_q > q >p$. Let $\lambda'$ be the image of $\lambda$ under the  transposition $(p,q)$, that is, with $j_p$ and $j_q$  exchanged. Then 
$s_{\lambda'} (\tb) > 0$. The coefficients of these two terms in 
\eqref{eq:sum}, $A_\lambda$ and $A_{\lambda'}$ differ by  two terms only.  We have  
$$
 A_{\lambda'} - A_\lambda = 
   (a_{p,j_q -p} a_{q,j_p -q} -  a_{p,j_p -p} a_{q,j_q -q})
     \prod_{i \ne  p,q} a_{i,j_i -i}.   
$$
We now show that $A_{\lambda'} > A_\lambda$, which will complete the
proof of the proposition. Recall the definition of $a_{k,j}$ in 
\eqref{eq:akj}. Let us denote $b_{k,j} = a_{k,j} /\binom{r-s}{j}$. 
Then it is easy to verify that 
$$
\frac{b_{p,j_q-p}}{b_{q,j_q-q}} = \frac{(2q-1)(2q-3) \ldots (2p+1)}
  {(2s-2q-1)(2s-2q-3)\ldots (2s-2p+1)}  =: B_{p,q},
$$
which is independent of $j_p$ and $j_q$. Consequently, we have 
\begin{align*}
 \frac{a_{p,j_q-p}}{a_{q,j_q-q}} & = B_{p,q} 
     \frac{\binom{r-s}{j_q-p}} {\binom{r-s}{j_q-q}} 
   = B_{p,q} \left(\frac{r-s+1}{j_q-q+1} -1\right) \ldots 
     \left(\frac{r-s+1}{j_q-p} -1\right) \\
  & >  B_{p,q} \left(\frac{r-s+1}{j_p-q+1} -1\right) \ldots 
     \left(\frac{r-s+1}{j_p-p} -1\right)  =\frac{a_{p,j_p-p}}{a_{q,j_p-q}}, 
\end{align*}
which implies that $A_{\lambda'} > A_\lambda$ and completes the proof.
\end{proof}

\begin{prop} \label{prop:2}
Let $m$ and $k$ be integers such that $1 \le k \le m$. Let $p_{2m-k-1}$ 
and $q_{k-1}$ be polynomials of degree $2m-k -1$ and degree $k-1$, 
respectively. If $t_1, \ldots, t_m$ be distinct numbers in $(0,1)$ and 
\begin{align} 
& p_{2m-k-1}^{{\rm even}}(t_i)+ q_{k-1}^{{\rm odd}}(t_i)
    (1-t_i^2)^{m-k} = 0, \label{eq:pq1}\\
& p_{2m-1-k}^{{\rm odd}}(t_i)+q_{k-1}^{{\rm even}}(t_i)(1-t_i^2)^{m-k}=0,
\label{eq:pq2} 
\end{align}
then $p_{2m-k-1}(t) \equiv 0$ and $q_{k-1}(t) \equiv 0$. 
\end{prop}

\begin{proof}
Depending on $k$ being even or odd, we need to consider the following four cases. 

\medskip\noindent
{\sc Case 1.} $k$ is even. Setting $r = m - (k+2)/2$ and $s = (k-2)/2$, 
 equation \eqref{eq:pq1} becomes 
$$
  p_r(t_i^2) + t_i \, q_s (t_i^2)(1-t_i^2)^{r-s} = 0, \qquad 1\le i \le r+s+2. 
$$
From Proposition \ref{prop:1} with $\epsilon =1$ and $t^{\pm 1} = t$,
it follows that $p_r(t) \equiv 0$ and $q_s(t) \equiv 0$. 

\medskip\noindent
{\sc Case 2.} $k$ is even. Setting $r = m - (k+2)/2$ and $s = (k-2)/2$, 
 equation \eqref{eq:pq2} becomes 
$$
  t_i \, p^*_r(t_i^2) +  q^*_s (t_i^2)(1-t_i^2)^{r-s} = 0, 
      \qquad 1\le i \le r+s+2. 
$$
Multiplying the equation by $t_i^{-1}$, we can use Proposition \ref{prop:1} 
with $\epsilon =1$ and $t^{\pm 1} = t^{-1}$ to conclude that 
$p_r(t) \equiv 0$ and $q_s(t) \equiv 0$. 

\medskip\noindent
{\sc Case 3.} $k$ is odd. Setting $r = m - (k+1)/2$ and $s = (k-1)/2$, 
 equation \eqref{eq:pq1} becomes 
$$
   p^*_r(t_i^2) + t_i q^*_{s-1} (t_i^2)(1-t_i^2)^{r-s} = 0, 
     \qquad 1\le i \le r+s+1. 
$$
From Proposition \ref{prop:1} with $\epsilon =0$ and $t^{\pm 1} = t$,
it follows that $p_r(t) \equiv 0$ and $q_s(t) \equiv 0$.  

\medskip\noindent
{\sc Case 4.} $k$ is odd. Setting $r = m - (k+2)/2$ and $s = (k-2)/2$, 
 equation \eqref{eq:pq2} becomes 
$$
  t_i \, p^*_r(t_i^2) +  q^*_{s-1} (t_i^2)(1-t_i^2)^{r-s} = 0, 
    \qquad 1\le i \le r+s+1. 
$$
Multiplying the equation by $t_i^{-1}$, we can use Proposition \ref{prop:1} 
with $\epsilon =1$ and $t^{\pm 1} = t$ to conclude that $p_r(t) \equiv 0$
and $q_s(t) \equiv 0$. 
\end{proof}

\subsection{Factorization method}
The following factorization theorem holds the key to our main result.

\begin{thm}\label{teo:fac}
Let $m$ and $s$ be positive integers  satisfying $m\le s\le 2m-1$. Denote 
$\lambda= s-m+1$. Let $\theta_1,\dots,\theta_{2\lambda}$ be distinct 
numbers in $(0,\pi)$ such that $\theta_{2\lambda+1-i}=\pi-\theta_i$ for 
$i=1,\dots,\lambda$. Denote 
$$
\wt X = \{ (\theta_i,\phi_{i,j}): \phi_{i,j} \in \Theta_{0,m},  
\,\, 1 \le i \le \lambda, \,\, \hbox{and} \,\, 
 \phi_{i,j} \in \Theta_{1,m}, \,\, \lambda+1 \le i \le 2 \lambda \}.
$$
If $T_s\in\Pi_{s}(\mathbb{S}^2)$ satisfies
\[
  \wt T_s(\theta_i,\phi_{i,j})=0,\qquad (\theta_i,\phi_{i,j})\in \wt X, 
\]
then there is a spherical polynomial $T^{\ast}_{s-2\lambda} \in
\Pi_{s-2\lambda}(\mathbb{S}^2)$ such that
\[
  T(x,y,z)=\prod_{i=1}^{2\lambda}(z-\cos\theta_i)ß,T_{s-2\lambda}^{\ast}(x,y,z).
\]
In particular, $T_{s-2\lambda}^{\ast}\equiv 0$ if $s\!=\!2m\!-\!1$. 
\end{thm}

\begin{proof} 
We start with the formula \eqref{eq:polynomial}, which becomes 
\begin{align*}
 \wt T_n(\theta,\phi)= a_0(\cos\theta) +\sum\limits_{k=1}^{s}
  \left[a_k(\cos\theta)\,(\sin\theta)^k\cos k\phi
         +b_k(\cos\theta)(\sin\theta)^k\sin k\phi \right], \notag
\end{align*}
where $a_k(\cdot)$ and $b_k(\cdot)$ are polynomials of degree $s-k$. 
For $i =1,2\ldots, \lambda$, we can follow the proof of Lemma \ref{lem:2.2} 
and Lemma \ref{lem:2.3} and distinguish the following three cases.

\medskip\noindent
{\sc Case 1.} For $0\le k\le 2m-s-1$,
\[
a_k(\cos\theta_i)(\sin\theta_i)^k=0,\quad b_k(\cos\theta_i)(\sin\theta_i)^k=0,
\quad i=1,\dots,2\lambda, 
\] 
setting $b_0\equiv 0$. Since $\theta_i\in (0,\pi)$, we have that 
$a_k(\cos\theta_i)=0 $ and $b_k(\cos\theta_i)=0$, for $i=1,\dots,2\lambda$. 
Recall that $a_k$ and $b_k$ are polynomials of degree $s-k\ge s-(2m-s-1)=
2\lambda-1$. Consequently, there exist polynomials $a_k^{\ast}$ and 
$b_k^{\ast}$, both of degree $s-k-2\lambda$, such that 
\[
a_k(t)=\prod_{i=1}^{2\lambda}(t-\cos\theta_i)a_k^{\ast}(t)
 \quad\hbox{and}\quad
b_k(t)=\prod_{i=1}^{2\lambda}(t-\cos\theta_i)b_k^{\ast}(t).
\] 
In the extreme case $k\!=\!2m\!-\!s\!-\!1$, we have $a_{2m-s-1}=b_{2m-s-1}=0$.

\medskip\noindent
{\sc Case 2.} For $k=m$,
\begin{align*}
& (\sin\theta_i)^m a_m(\cos\theta_i)=0\quad \hbox{and}\quad
      i=1,\dots,\lambda,\\
& (\sin\theta_i)^m b_m(\cos\theta_i)=0\quad \hbox{and}\quad
      i=\lambda+1,\dots,2\lambda.
\end{align*}
Since $\theta_i\in (0,\pi)$ and both $a_m$ and $b_m$ are polynomials of 
degree $s-m=\lambda-1$  vanishing at at least $\lambda$ points, 
they have to be identically zero.

\medskip\noindent
{\sc Case 3.} For $2m-s\le k\le m-1$, we end up with equations similar to  \eqref{eq:2.6}:
\begin{align} \label{eq:final} 
\begin{split}
& p_{s-k}^{{\rm even}}(t_i)+ q_{2\lambda - (s-k)-2}^{{\rm odd}}(t_i)
    (1-t^2)^{s-\lambda-k+1} = 0, \\
& p_{s-k}^{{\rm odd}}(t_i) + 
   q_{2 \lambda -(s-k)-2}^{{\rm even}}(t_i)(1-t^2)^{s-\lambda-k+1}=0, 
\end{split}
\end{align}
for $i = 1,2,\ldots, \lambda$, where either $p_{s-k}(t) = a_k(t)$ and
$q_{2\lambda-(s-k)-2}(t) = a_{2m-k}(t)$, or $p_{s-k}(t) = b_k(t)$ 
and $q_{2\lambda-(s-k)-2}(t)= - b_{2m-k}(t)$. In deriving the above
equations we have used several times the identity $\lambda = s-m+1$. Recall that 
$a_k$ and $b_k$ are polynomials of degree $s-k$; the subscript of the 
polynomials $p_{2\lambda-k}$ and $q_{2\lambda-(s-k)-2}$ again indicate 
their degree. 

It is easy to see that the system of equations \eqref{eq:final} is 
exactly the one being  studied in the previous subsection, namely \eqref{eq:pq1} and \eqref{eq:pq2}. Hence, 
using Proposition \ref{prop:2} we conclude that 
\begin{align*}
& a_k(t) \equiv 0 \quad\hbox{and}\quad 
   a_{2m-k}(t) \equiv 0, \quad k=2m-s,\dots,m-1, \\
& b_k(t) \equiv 0 \quad\hbox{and}\quad
    b_{2m-k}(t) \equiv 0, \quad k=2m-s,\dots,m-1.
\end{align*}
Together, these three cases show that we have the factorization
\begin{align*}
\wt T_s(\theta,\phi) & = \prod_{i=1}^{2\lambda} (\cos\theta-\cos\theta_i)\\
&\times\left(a_0^{\ast}(\cos\theta)+\sum_{k=1}^{2m-s-2}(a_k^{\ast}(\cos\theta)\cos k\phi+b_k^{\ast}(\cos\theta)\sin k\phi)\right),
\end{align*}
which completes the proof.
\end{proof}

Using  factorization repeatedly, we can obtain a complete factorization 
of a polynomial of degree $2m-1$ in $\Pi_{2m-1}(\mathbb{S}^2)$.

\begin{thm}\label{teo:fac1}
Let $n$ be an odd positive integer, $\sigma\in\mathbb{N}$, and 
$\lambda_1,\dots,\lambda_{\sigma}$ be positive integers. Define 
$n_k=n_{k-1}-2\lambda_k$, for $1 \le k \le \sigma$, with $n_0 = n$. 
Assume that $n_k \ge 0$ for $1 \le k \le \sigma -1$. If $T_n 
\in \Pi_n(\mathbb{S}^2)$ satisfies 
$$
  \wt T_n(\theta_{i,k},\phi_{i,j,k}) =0, \quad 1 \le i \le 2 \lambda_k, \,\,
 0 \le j \le 2 (n_{k-1}-\lambda_k+1)-1, \,\, 1 \le k \le \sigma, 
$$
where $\theta_{i,k}$ are pairwise
 distinct angles in $(0, \pi)$ with $\theta_{i,2\lambda_k+1-l}\!=\!
\pi-\theta_{i,l}$, $l=1,\dots,\lambda_k$, $\phi_{i,j,k} \in 
\Theta_{0,n_{k-1}-\lambda_k+1}$ for $1 \le i \le \lambda_k$ and 
$\phi_{i,j,k} \in  \Theta_{1,n_{k-1}-\lambda_k+1}$ for 
$\lambda_k+1 \le i \le 2 \lambda_k$, then there exists a polynomial 
$T_{n_{\sigma}}^* \in \Pi_{n_\sigma}(\mathbb{S}^2)$ such that 
$$
 T_n(x,y,z) = \prod_{k=1}^\sigma \prod_{i=1}^{2\lambda_k} 
    (z-\cos \theta_{i,k})\, T_{n_{\sigma}}^*(x,y,z).
$$
In particular, $T_n(x,y,z) \equiv 0$ if $n_\sigma < 0$.
\end{thm}

\begin{proof}
We apply the factorization result in Theorem \ref{teo:fac} repeatedly with
$s = n_{k-1}$, $m = n_{k-1} - \lambda_k+1$ and $\lambda = \lambda_k$ for 
$k = 1, 2, \ldots, \sigma$.  
\end{proof}  

Just as in  the case of an odd number of points on each latitude (cf. \cite{X1}), the interpolation nodes in the above theorem are located on latitudes split up into $\sigma $ groups $\{\mathbb{S}^2(z_{i,k}): 1 \le i \le 
2\lambda_k\}$, $1 \le k \le \sigma$, $z_{i,k} = \cos \theta_{i,k}$, and 
$z_{2\lambda_k+1-i,k}=-z_{i,k}$. Latitudes in different groups contain 
 a different number of nodes. More precisely, each of the latitudes in the 
$k$-th group, $\mathbb{S}^2(z_{1,k}), \mathbb{S}^2(z_{2,k}), \ldots, 
\mathbb{S}^2(z_{2\lambda_k,k})$, contains an even number of 
$2 (n_{k-1}- \lambda_k+1) $ equidistant points and the points lie on 
symmetric latitudes. In other words,  points on $\mathbb{S}^2(z_{2\lambda_k+1-i,k})$ 
and $\mathbb{S}^2(z_{i,k})$, $i=1,\dots,\lambda_k$, differ by a rotation of an angle of  $\pi/(2(n_{k-1}-\lambda_k+1))$.

\section{Interpolation on the sphere }
\setcounter{equation}{0}

Our main result on interpolation follows from the factorization 
Theorem \ref{teo:fac1}. The following formula can be used to verify that 
the number of interpolation conditions matches the dimension of the polynomial
space: 
\begin{equation}\label{eq:dimension}
 \Pi_s(\mathbb{S}^2) = \dim \Pi_{s-2\lambda }(\mathbb{S}^2) + 
    2\lambda\,(2s-2\lambda+2).
\end{equation}

\begin{thm} \label{thm:3.1}
Let $n$ be an odd natural number and let $\lambda_1, \ldots, \lambda_\sigma$ $(\sigma\in\mathbb{N})$ be positive integers, such that  
\begin{align} \label{eq:3.2}
   \lambda_1 + \ldots + \lambda_\sigma = \frac{n+1}{2}.
\end{align}
Define $n_k = n_{k-1} - 2 \lambda_k $, for $1 \le k \le \sigma -1$, with
$n_0= n$. Let 
$$
\wt  X = \{(\theta_{i,k},\phi_{i,j,k}): 1 \le i \le 2 \lambda_k, \,\,
 0 \le j \le 2 (n_{k-1}-\lambda_k+1)-1, \,\, 1 \le k \le \sigma \}, 
$$
where $\theta_{i,k}$, $1 \le j \le 2 \lambda_k$ and $1 \le k \le \sigma$, are
distinct numbers in $(0, \pi)$ with $\theta_{2\lambda_k+1-i,k}=
\pi-\theta_{i,k} $ $(i=1,\dots,\lambda_k)$, $\phi_{i,j,k} \in 
\Theta_{0,n_{k-1}-\lambda_k+1}$, for $1 \le i \le \lambda_k$, and 
$\phi_{i,j,k} \in  \Theta_{1,n_{k-1}-\lambda_k+1}$, for 
$\lambda_k+1 \le i \le 2 \lambda_k$. Then the set $X$ solves the 
interpolation problem in $\Pi_n(\mathbb{S}^2)$.
\end{thm}

\begin{proof}
First, we verify that the dimension of $\Pi_n(\mathbb{S}^2)$ matches 
the number of  interpolation conditions. Let $|X|$ denote the
number of points in $X$. It follows from  equation \eqref{eq:dimension} 
that 
\begin{align*}
|X| = & \sum_{k=1}^\sigma 2 \lambda_k\,(2 n_{k-1} - 2 \lambda_k+2)\\
= & \sum_{k=1}^{\sigma}(\dim \Pi_{n_{k-1}}(\mathbb{S}^2) -
 \dim \Pi_{n_{k-1}-2\lambda_k}(\mathbb{S}^2))=(n+1)^2
=\dim \Pi_n(\mathbb{S}^2).
\end{align*}
 Thus, it is sufficient to show that if 
$T_n \in \Pi_n(\mathbb{S}^2)$ vanishes on $X$, then $T_n(x,y,z) \equiv 0$. 
Under the condition \eqref{eq:3.2}, it follows that 
$$
n_\sigma = n_{\sigma -1} - 2\lambda_\sigma  = 
 n_{\sigma-2}-2\lambda_{\sigma-1}-2\lambda_{\sigma}= \cdots 
 =n - 2\sum_{k=1}^{\sigma-1}\lambda_k=-1  < 0.   
$$
Hence, the factorization in Theorem \ref{teo:fac1} shows that 
$T(x,y,z) \equiv 0$.
\end{proof} 

For a fixed $n =2m-1$, this theorem contains a number of different 
interpolation processes. In fact, for each positive odd integer $n$, the 
number of point sets $X$ which can be deduced from  Theorem \ref{thm:3.1} depends on the 
partition number of $(n+1)/2$. Every solution of  equation 
\eqref{eq:3.2} leads to a set of  points  defining a poised  
interpolation problem in $\Pi_n(\mathbb{S}^2)$. The number of solutions
of such an equation grows exponentially as $n$ goes to infinity. Moreover, 
the order of $\lambda_1, \ldots, \lambda_\sigma$ matters; i.e. , different 
permutations of a solution $\lambda_1, \ldots, \lambda_\sigma$ of  equation
\eqref{eq:3.2} give different sets of interpolation points. 

Among the solutions of equation  \eqref{eq:3.2}, one extreme case is
$\sigma = 1$, for which the equation has only one solution $\lambda_1 = 
(n+1)/2$. In this case, the interpolation points are located 
on $n+1$ symmetric latitudes $\mathbb{S}^2(z_1)$, $\mathbb{S}^2(z_2), \ldots, 
\mathbb{S}^2(z_{n+1})$, each of them  containing  $n+1$ equidistant points. 
This case has already appeared in \cite[Theorem 2.5]{F}.

The other extreme case is $\sigma=(n+1)/2$ and $\lambda_1=\cdots=
\lambda_{(n+1)/2}=1$. There,  $n_k = n-2k$ and  we have $(n+1)/2$ 
groups of two symmetric  latitudes, where the ones in the  $k$th group 
contain $2(n-2k+2)$ equidistant nodes. The points on  a latitude 
 are rotated  by an 
angle of $\pi/(2(n-2k+2))$ with respect to the points on the corresponding symmetric latitude of the same group.  

\begin{exam} To illustrate the power of the factorization method, 
we present the possible point distributions for $n=3,5$ and $7$.
\begin{itemize}
\item $n=3$
\begin{enumerate}
\item $\sigma=1,  \lambda_1=2$: $4$ latitudes each  with $4$ points,
\item $\sigma=2,\lambda_1=\lambda_2=1$: $2$ latitudes with $6$ points and
 $2$ latitudes with $2$ points. 
\end{enumerate}
\item $n=5$
\begin{enumerate}
\item $\sigma=1,\lambda_1=3$: $6$ latitudes, each with $6$ points.
\item $\sigma=2,\lambda_1+\lambda_2=3$ has two solutions.
\begin{enumerate}
\item $\lambda_1=1,\lambda_2=2$: $4$ latitudes with $8$ points and $2$ 
latitudes with $2$ points.
\item $\lambda_1=2,\lambda_2=1$: $2$ latitudes with $10$ points and
 $4$ latitudes with $4$ points.
\end{enumerate}
\item $\sigma=3, \lambda_1=\lambda_2=\lambda_3=1$: $2$ latitudes with 
$10$ points, $2$ latitudes with $6$ points and $2$ latitudes with $2$ points.
\end{enumerate}
\item $n=7$
\begin{enumerate}
\item $\sigma=1,\lambda_1=4$: $8$ latitudes with $8$ points.
\item $\sigma=2, \lambda_1+\lambda_2=4$ has three solutions,
\begin{enumerate}
\item $\lambda_1=2,\lambda_2=2$:  $4$ latitudes with $12$ points and 
$4$ latitudes with $4$ points;
\item $\lambda_1=1,\lambda_2=3$:  
$2$ latitudes with $14$ points and $6$ latitudes with $6$ points;
\item $\lambda_1=3,\lambda_2=1$:  $6$ latitudes with $10$ points 
and $2$ latitudes with $2$ points.
\end{enumerate} 
\item $\sigma=3, \lambda_1+\lambda_2+\lambda_3=4$ has three solutions,
\begin{enumerate}
\item $\lambda_1=1,\lambda_2=1,\lambda_3=2$:  $2$ latitudes with $14$ points,
 $2$ latitudes with $10$ points and  $4$ latitudes with $4$ points.
\item $\lambda_1=1,\lambda_2=2,\lambda_3=1$:  $2$ latitudes with $14$ points,
 $4$ latitudes with $8$ points and  $2$ latitudes with $2$ points;
\item $\lambda_1=2,\lambda_2=1,\lambda_3=1$:  $4$ latitudes with $12$ points,
 $2$ latitudes with $6$ points and  $2$ latitudes with $2$ points.
\end{enumerate} 
\item $\sigma=4,\lambda_1=\lambda_2=\lambda_3=\lambda_4=1$: $2$ latitudes 
with $14$ points, $2$ latitudes with $10$ points, $2$ latitudes with $6$ 
points and $2$ latitudes with $4$ points.
\end{enumerate}
\end{itemize}
\end{exam}

It is well-known that interpolating polynomials can be used to construct
cubature formulas on the unit sphere (cf. \cite{S}). In fact,  integrating   the interpolation polynomial in $\Pi_n(\mathbb{S}^2)$ yields a cubature formula on the sphere 
which is exact for spherical polynomials of degree $n$. Among the point
sets in Theorem \ref{thm:3.1}, the case where the points are distributed 
on $2m$ symmetric latitudes, with each latitude containing $2m$ equidistant 
points, is of particular interest. In this case, the cubature formula 
is simple and can be explicitly given. 

\begin{prop}
Let $m$ be a positive integer. Let $\theta_1, \ldots, \theta_{2m}$ be 
pairwise distinct numbers in $(0,\pi)$ with $\theta_{2m+1-i}=\pi-\theta_i$, $i=1,\dots,m$, and  $\alpha\in\lbrace 0,1\rbrace$. Then 
for all $T_{2m-1} \in \Pi_{2m-1}(\mathbb{S}^2)$, 
\begin{align*}
 \int_{\mathbb{S}^2} T_{2m-1}(\xi)\, d\omega(\xi)  =
  \frac{\pi}{m} \sum_{i=1}^{m} \lambda_i
   \sum_{j=0}^{2m-1} \wt T_{2m-1}(\theta_i,\phi_{j}^{0}) 
 + \frac{\pi}{m} \sum_{i=m+1}^{2m} \lambda_i
   \sum_{j=0}^{2m-1} \wt T_{2m-1}(\theta_i,\phi_{j}^{1})
\end{align*}
where $\phi_{j}^{\alpha} = (2j+\alpha)\pi/2m$, and $\lambda_i$ is given by 
$$
\lambda_i =  \int_{-1}^1 \prod_{k=1,k\ne i}^{2m} 
 \frac{t-\cos \theta_k}{\cos \theta_i - \cos \theta_k} dt,\qquad i=1,\dots,2m.
$$
\end{prop} 

\begin{proof}
Let the interpolation polynomial $T_{2m-1}$ be of the form 
\eqref{eq:polynomial}. We use the quadrature formula 
\begin{equation} \label{eq:trigquad}
  \frac{1}{2 \pi} \int_0^{2\pi} p(t) dt = \frac{1}{2m} 
  \sum_{j=0}^{2m-1} p(\phi_j^\alpha), 
\end{equation}
which is known to hold for every trigonometric polynomial of degree at most
$m$ (see, for example, \cite[Vol.2, p. 8]{Z}). Using  formula 
\eqref{eq:trigquad} and the interpolation property of $T_{2m-1}$, it 
follows that 
$$
a_0(\cos \theta_i) = \frac{1}{2\pi}\int_0^{2\pi} \wt T_{2m-1}(\theta_i,\phi)\, 
 d\phi  
= \frac{1}{2m} \sum_{j=0}^{2m-1} \wt  T_{2m-1}(\theta_i,\phi_{j}^{\alpha}), 
$$
for every fixed $\theta_i$, $1\le i \le 2m$. 
Consequently, $a_0$, which is a polynomial of degree $2m-1$ in one variable,
is uniquely determined by these $2m$ interpolation conditions. It follows that 
$$
a_0(t) = \sum_{i=1}^{m} \left(\frac{1}{2m}\sum_{j=0}^{2m-1} 
        \wt T_{2m-1}(\theta_i,\phi_{j}^{0}) \right) \ell_i(t)
   + \sum_{i=m+1}^{2m} \left(\frac{1}{2m}\sum_{j=0}^{2m-1} 
        \wt T_{2m-1}(\theta_i,\phi_{j}^{1}) \right) \ell_i(t),  
$$
where $\ell_i(t) = \prod_{k=1,k\ne i}^{2m} 
            (t-\cos \theta_k)/(\cos \theta_i - \cos \theta_k)$.
Using the change of variables 
$$
 \int_{\mathbb{S}^2} T_{2m-1}(\xi) d\omega(\xi) = \int_0^\pi \int_0^{2\pi} \wt T_{2m-1}(\theta,\phi)
     \sin \theta\, d\phi\, d \theta,
$$
the integral of $T_{2m-1}$ over the surface of the sphere is equal to 
$$
\int_{\mathbb{S}^2} T_{2m-1}(\xi) d \omega(\xi) = 
  2 \pi \int_0^\pi a_0(\cos \theta) 
        \sin \theta d \theta = 2 \pi \int_{-1}^1 a_0(t) dt. 
$$
The stated formula follows from the formula for $a_0(t)$ given above.
\end{proof}

In particular, this result shows that the cubature formula is nonnegative, 
if $\cos \theta_i$ are chosen so that $\lambda_i$ are nonnegative. This holds,
for example, if $\cos \theta_i$ are the zeros of the Legendre polynomial 
$P_{2m}$ of degree $2m$, or the  zeros of a quasi Legendre orthogonal 
polynomial $P_{2m} + \alpha P_{2m-1}$ with mild conditions imposed on 
$\alpha \in \RR$ (see, for example, \cite{X94}). In \cite{PR},
the positivity of the cubature in this case has been  proved by working directly with the 
interpolation matrix.

\end{document}